\documentclass{amsart}

\numberwithin{equation}{section}
\def\textrm#1{\text{\textnormal{#1}}}

\theoremstyle{plain}
\newtheorem{theorem}{Theorem}[section]
\newtheorem{prop}[theorem]{Proposition}

\newtheorem{lemma}[theorem]{Lemma}

\theoremstyle{definition}

\theoremstyle{remark}
\newtheorem{remark}[theorem]{Remark}

 \def\today{\ifcase\month\or
  January\or February\or March\or April\or May\or June\or
  July\or August\or September\or October\or November\or December\fi
  \space\number\day, \number\year}

\begin{document}
\title[Division algebras]
{Division algebras that ramify only on a plane quartic curve}
\author[Kunyavski\u\i , Rowen, Tikhonov, Yanchevski\u\i ]
{Boris \`E. Kunyavski\u\i , Louis H. Rowen, Sergey V. Tikhonov,
and \\ Vyacheslav I. Yanchevski\u\i }

\thanks{This research was supported by
the Israel Science Foundation founded by the Israel Academy of
Sciences and Humanities --- Center of Excellence Program and by
RTN Network HPRN-CT-2002-00287.}
\thanks{The first author was partially supported by the Ministry of
Absorption (Israel) and the Minerva Foundation through the Emmy
Noether Research Institute of Mathematics.}
\thanks{The third and the fourth authors were partially
supported by the Fundamental Research Foundation of Belarus.}

\date{\today}
\address{Kunyavski\u\i , Rowen:
Department of Mathematics, Bar-Ilan University,
52900 Ramat Gan, ISRAEL} \email{kunyav@macs.biu.ac.il ,
rowen@macs.biu.ac.il}
\address{Tikhonov, Yanchevski\u\i : Institute of Mathematics
of the National Academy of Sciences of Belarus, ul. Surganova 11,
220072 Minsk, BELARUS} \email{tsv@im.bas-net.by ,
yanch@im.bas-net.by}

\begin{abstract}
Let $k$ be an algebraically closed field of characteristic 0. We
prove that any division algebra over $k(x,y)$ whose ramification
locus lies on a quartic curve is cyclic.
\end{abstract}

\maketitle

\def\dd{{\partial}}

\def\toeq{{@>\sim>>}}
\def\into{{\hookrightarrow}}

\def\emptyset{{\varnothing}}

\def\alp{{\alpha}}  \def\bet{{\beta}} \def\gam{{\gamma}}
 \def\del{{\delta}}
\def\eps{{\varepsilon}}
\def\kap{{\kappa}}                   \def\Chi{\text{X}}
\def\lam{{\lambda}}
 \def\sig{{\sigma}}  \def\vphi{{\varphi}} \def\om{{\omega}}
\def\Gam{{\Gamma}}  \def\Del{{\Delta}}  \def\Sig{{\Sigma}}
\def\ups{{\upsilon}}


\def\A{{\mathbb A}}
\def\F{{\mathbb F}}
\def\Q{{\mathbb{Q}}}
\def\CC{{\mathbb{C}}}
\def\PP{{\mathbb P}}
\def\R{{\mathbb R}}
\def\Z{{\mathbb Z}}
\def\X{{\mathbb X}}

\def\Fb{{\overline F}}
\def\Hb{{\overline H}}
\def\Kb{{\overline K}}
\def\Lb{{\overline L}}
\def\Yb{{\overline Y}}
\def\Xb{{\overline X}}
\def\Tb{{\overline T}}
\def\Bb{{\overline B}}
\def\Gb{{\overline G}}
\def\Vb{{\overline V}}

\def\kb{{\bar k}}
\def\xb{{\bar x}}

\def\Th{{\hat T}}
\def\Bh{{\hat B}}
\def\Gh{{\hat G}}

\def\Xt{{\tilde X}}
\def\Gt{{\tilde G}}

\def\gg{{\mathfrak g}}
\def\gm{{\mathfrak m}}
\def\gp{{\mathfrak p}}
\def\gq{{\mathfrak q}}

\def\min{^{-1}}

\def\textrm#1{\text{\rm #1}}

\def\char{\textrm{char}}
\def\cor{\textrm{cor}}
\def\ord{\textrm{ord}}
\def\ram{\textrm{ram}}
\def\deg{\textrm{deg}}
\def\exp{\textrm{exp}}
\def\Gal{\textrm{Gal}}
\def\Spec{\textrm{Spec}}
\def\Proj{\textrm{Proj}}
\def\Perm{\textrm{Perm}}
\def\coker{\textrm{coker\,}}
\def\Hom{\textrm{Hom}}
\def\im{\textrm{im\,}}
\def\ind{\textrm{ind}}
\def\int{\textrm{int}}
\def\inv{\textrm{inv}}

\def\tors{_{\textrm{tors}}}      \def\tor{^{\textrm{tor}}}
\def\red{^{\textrm{red}}}         \def\nt{^{\textrm{ssu}}}
\def\sc{^{\textrm{sc}}}
\def\sss{^{\textrm{ss}}}          \def\uu{^{\textrm{u}}}
\def\ad{^{\textrm{ad}}}           \def\mm{^{\textrm{m}}}
\def\tm{^\times}                  \def\mult{^{\textrm{mult}}}
\def\tt{^{\textrm{t}}}
\def\uss{^{\textrm{ssu}}}         \def\ssu{^{\textrm{ssu}}}
\def\cf{^{\textrm{cf}}}
\def\ab{_{\textrm{ab}}}

\def\et{_{\textrm{\'et}}}
\def\nr{_{\textrm{nr}}}

\def\op{^{\textrm{op}}}

\def\til{\;\widetilde{}\;}

\def\eqf{{\sim_{\mathcal{F}}}}

\def\emptyset{{\varnothing}}
\def\div{{\textrm{div}}}

\def\cA{{\mathcal A}}
\def\cB{{\mathcal B}}
\def\cC{{\mathcal C}}
\def\cD{{\mathcal D}}
\def\cV{{\mathcal V}}

\font\cyr=wncyr10 scaled \magstep1%
\def\Bcyr{\text{\cyr B}}
\def\Sh{\text{\cyr Sh}}
\def\Ch{\text{\cyr Ch}}

\def\lpsi{{{}_\psi}}
\def\bks{{\backslash}}

\def\Br{\textrm{Br}}
\def\Pic{\textrm{Pic}}
\def\Bt{{{}_2\textrm{Br}}}
\def\Bn{{{}_n\textrm{Br}}}


\section{Main result} \label{sec:intro}

Let $k$ be an algebraically closed field of characteristic 0,
$L=k(x,y)$ its purely transcendental extension of degree 2. We are
interested in division $L$-algebras with ramification only at a
quartic curve. Our main result (Theorem \ref{th:main}), which grew
out from our earlier paper \cite{KRTY}, shows that any such
algebra is a symbol algebra whose index equals the exponent. This
is a generalization of earlier results of Ford, Van den Bergh and
Jacob on cyclicity of $k(x,y)$-algebras ramified either along a
cubic curve, or a hyperelliptic curve, or a special quartic curve
(see \cite{Fo}, \cite{Fo2}, \cite{Fo3}, \cite{Ja}, \cite{vdB}). In
contrast to the above cited papers (and to the paper of de Jong
\cite{dJ} who proved the equality of index and the exponent for
any central simple algebra defined over an extension of $k$ of
transcendence degree 2, regardless of its ramification locus), our
methods are quite elementary and describe the cyclic algebra
explicitly in terms of the ramification. This might be a
justification for publishing this result, in spite of the fact
that the long-standing problem of cyclicity of $k(x,y)$-algebras
seems to be close to a general solution by M.~Ojanguren and
R.~Parimala (private communication of Parimala to one of the
authors). Roughly, our approach consists in regarding $L$ as the
function field $l({\mathbb P}^1_l)$ of the projective {\it line}
over $l=k(x)$ (rather than the function field $k({\mathbb P}^2_k)$
of the projective {\it plane} over $k$). As observed by the
anonymous referee of \cite{KRTY}, this allows one to apply the
full strength of the Faddeev theory of algebras over a projective
line (and, in particular, to employ classification results of
\cite{KRTY} for algebras with small ramification locus).

\bigskip

We fix an integer $n > 1$. Let $\Br(F)$ denote the Brauer group of
a field $F$ and $_n\Br(F)$ its $n$-torsion subgroup.
Assume $F$ contains a primitive $n$th root of unity $\rho$ whose
notation is fixed. For any finite dimensional $F$-central simple
algebra $\cA$, we use $[\cA]$ to denote its class in $\Br(F)$. The
{\it exponent} $\exp(\cA)$ of $\cA$ is the order of $[\cA]$ in
$\Br(F)$, and the {\it index} $\ind(\cA)$ of $\cA$ is the degree
of the central division algebra $\cD$ such that $\cA\cong
M_r(\cD)$. If $a,b\in F^*$, we denote by $(a,b)_n$ the cyclic
algebra generated over $F$ by $\alp ,\beta$ such that $\alp ^n=a$,
$\beta ^n=b$, $\alp\beta =\rho\beta\alp$. (We often drop the
subscript $n$ if this does not lead to any confusion.)

Let $R$ be a discrete valuation ring with quotient field $K$ and
residue field $l$. We assume that both $l$ and $K$ are of
characteristic zero. Then it is classical \cite{AB} that there is
an exact sequence

$$
0 \rightarrow \Br(R) \rightarrow \Br(K)
\stackrel{\ram}{\longrightarrow} H^1(l,\Q/\Z) \rightarrow 0.
$$

The ramification map, denoted by $\ram$, is described in \cite{AB}
(see also \cite[Ch.10]{Sa99}).

Let now $k$ be an algebraically closed field of characteristic
zero, and let $L=k(\PP_k^2)$ be the function field of the
projective plane over $k$.

Let $C\subset \PP_k^2$ be an irreducible curve. The local ring
$\mathcal{O}_C$ is a discrete valuation ring, and we denote by
$v_C\colon L^* \rightarrow \Z$ the corresponding discrete
valuation.

We view $\PP_k^2 = \Proj (k[x,y,z])$ as $\Spec(k[x,y])\cup
\{L_\infty\}$ where $L_\infty$ is the infinite line corresponding
to $z = 0$. (Abusing notation, we write $k[x,y]$ for what really
should be $k[(x/z),(y/z)]$.) Except for $L_\infty$, each $C$
corresponds to an irreducible polynomial in either $k(x)[y]$ or
$k(y)[x]$, and the valuation $v_C$ comes from the localization at
that polynomial. At $L_\infty$, we must use the $1/x$-adic
valuation on $k(x/y)[x]$ to compute $v_{L_\infty}$.

Then there is the map corresponding to $\mathcal{O}_C$
$$
\ram_C\colon\Br(L) \rightarrow H^1(k(C),\Q/\Z),
$$
where $k(C)$ is the function field of $C$.

Note that $_n H^1(k(C),\Q/\Z)\cong k(C)^*/k(C)^{*n}$. By the
``tame symbol'' formula (see, for example, \cite[1.3 and 2.3]{Bl}),
the ramification of $[(a,b)_n]$ is
defined by the $n$th root of the residue of
$$(-1)^{v_C(a)v_C(b)}a^{v_C(b)}/b^{v_C(a)}.$$

According to \cite{AM}, there is a monomorphism $\Br(L) \to\oplus
H^1(k(C),\Q/\Z)$, where $C$ runs over all irreducible curves in
$\PP^2_k$; this map will also be denoted by $\ram$. For $[\cA]\in
\Br(L)$ we write
$$
R=\bigcup_C\{C\subset \PP_k^2 \;|\;\ram_C([\cA])\neq 0\}
$$
and call $R$ the {\it ramification divisor} of $\cA$.

Let $l$ be a field of characteristic zero. According to
{\textrm{\cite{Fa}, \cite[Ch. II, App., \S 5]{Se2}}}, there is an
exact sequence
\begin{equation} \label{eq:ramFad}
0\to {}_n\Br (l)\to {}_n\Br (l(\PP_l^1))\stackrel{\oplus\dd
_M}{\longrightarrow } \bigoplus _{M\in\PP ^1_l}H^1(l(M),\Z /n )
\stackrel{\oplus\cor_M}{\longrightarrow} H^1(l,\Z /n)\to 0.
\end{equation}
Here $M$ runs over the closed points of $\PP ^1_l$, and
$\oplus\cor_M$ is the sum of corestriction homomorphisms
$\cor_M\colon H^1(l(M),\Z /n )\longrightarrow H^1(l,\Z /n)$.
Moreover, using the isomorphism $H^1(l(M),\Z /n )\cong
l(M)^*/(l(M)^*)^n$ one can compute $\cor_M$ by the formula
$$\cor_M(a(l(M)^*)^n)=N_{l(M)/l}(a)(l^*)^n.$$
Sequence (\ref{eq:ramFad}) is sometimes referred to as the Faddeev
reciprocity law; indeed, it says that the sum of corestrictions of
the residues of a given algebra is zero.

Two $l(t)$-algebras $\cA$ and $\cB$ are said to be Faddeev
equivalent (our notation for that is $\cA\eqf\cB$) if there exists
an $l$-algebra $\cC$ such that $\cA$ is Brauer equivalent to
$\cB\otimes (\cC\otimes _ll(t))$. By the reciprocity law, two
algebras having the same residues in all but one closed rational
point of ramification must have the same residue in that point as
well and are thus Faddeev equivalent. In particular, two algebras
having the same residues in all finite points are Faddeev
equivalent.

\begin{remark} \label{rem:Tsen}
In case $\Br (l)$ is trivial, Faddeev equivalent algebras are
Brauer equivalent. This is the case when $l=k(x)$ for $k$
algebraically closed, by Tsen's theorem.
\end{remark}

There is a one-to-one correspondence between closed points of $\PP
^1_l$ and discrete valuations of $l(\PP_l^1)=l(t)$ trivial on $l$.
Recall the structure of discrete valuations of $l(t)$ trivial on
$l$. Any such valuation is of the following form. If it
corresponds to a finite closed point $M$ of $\PP ^1_l$, then there
exists an irreducible monic polynomial $f(t)\in l[t]$ such that
the valuation $v_f$ with $f(t)$ as a uniformizer coincides with
the valuation corresponding to $M$. The valuation corresponding to
the infinite point is $v_\infty$ with $t^{-1}$ as a uniformizer.

Let $f(t),g(t)\in l[t]$, $f$ irreducible and not dividing $g$.
Then the ramification of $(g(t),f(t))_n$ at $f$ is
$l(\theta)(\sqrt[n]{g(\theta)})/l(\theta)$ (here and throughout
below we identify elements of $l(\theta )^*/l(\theta )^{*n}$ with
degree $n$ cyclic extensions of $l(\theta )$).

\bigskip

Let $F(x,y)\in k[x,y]$ be an irreducible polynomial,
$F(x,y)=\sum_{i=0}^n g_i(x)y^i$, where $g_i(x)\in k[x]$,
$g_n(x)\ne 0$. Let $C$ be the curve defined by the equation
$F(x,y)=0$. Assume that $n>0$. As mentioned above, the field
$k(\PP_k^2)=k(x,y)$ can be interpreted as the function field
$k(x)(\PP_{k(x)}^1)=k(x)(y)$ of the line over the field $k(x)$.

Let $f(y)=F(x,y)/g_n(x)\in k(x)[y]$. Then $f(y)$ is a monic
irreducible polynomial over $k(x)$, and the valuation of $k(x,y)$
corresponding to the curve $C$ coincides with the valuation of
$k(x)(y)$ corresponding to the closed point defined by the
polynomial $f(y)$.

Thus, a central simple algebra $\cA/k(\PP_k^2)$ ramifies at the
valuation corresponding to the curve $C$ if and only if
$\cA/k(x)(\PP_{k(x)}^1)$ ramifies at a closed point defined by the
polynomial $f(y)$.

Note that if $n=0$ (i.e. $F(x,y)=x+c$, where $c\in k$), then the
valuation corresponding to the curve $C$ is not trivial on $k(x)$.
Hence in this case there is no valuation of $k(x)(y)$
corresponding to the valuation of $k(x,y)$ defined by $C$.
Moreover, an algebra $\cA$, viewed as an algebra over
$k(x)({\mathbb P}_{k(x)}^1)=k(x)(y)$, may ramify at the infinite
point even if $\cA$, viewed as an algebra over $k({\mathbb
P}_k^2)$, does not ramify at the infinite line $L_{\infty}$ given
by $z=0$. For example, consider the quaternion algebra $\cB
=(xy,x^2+1)$ \cite{KRTY}. It does not ramify at the infinite line
because the valuation $v_{L_{\infty }}$ (corresponding to the
infinite line) of both $xy$ and $x^2+1$ equals $-2$. Indeed,
$v_{L_{\infty }}$ coincides with the valuation $v_{\infty }$ of
$k(x/y)(x)$ corresponding to the infinite point. Hence we have
$$
v_{L_{\infty }}(xy)=v_{\infty}(xy)=v_{\infty}(x^2y/x)=-2
$$
since $y/x$ is ``constant''. Similarly, $ v_{L_{\infty
}}(x^2+1)=v_{\infty}(x^2+1)=-2$.

On the other hand, as an algebra over $k(x)({\mathbb
P}_{k(x)}^1)$, $\cB$ ramifies at the infinite point because
$x^2+1\notin k(x)^2$.

\bigskip

We are now ready to formulate our main result.

\begin{theorem} \label{th:main}
Let $k$ be an algebraically closed field of characteristic zero,
and $\cD$ a division algebra over $L=k(x,y)=k(\PP_k^2)$ with
ramification only along a quartic curve $C\subset \PP_k^2$. Then
$\cD$ is a cyclic algebra, and $\exp(\cD)=\ind(\cD)$.
\end{theorem}

The general strategy of the proof is as follows. As mentioned
above, we view $L$ as the function field $l(\PP^1_l)$, where
$l=k(x)$. Then we can use the following key result describing
central simple $k(x)$-algebras with small ramification locus up to
Faddeev equivalence:

\begin{prop} \label{prop:key}
Let $K$ be an arbitrary field of characteristic zero, containing
the $n$th roots of unity, $F=K(t)$, $\cA$ a central simple
$F$-algebra of exponent $n$. Suppose that $\cA$ has ramification
either

(i) at most at three linear points, or

(ii) at a point corresponding to a quadratic irreducible
polynomial, and perhaps at infinity, or

(iii) at a point corresponding to a cubic irreducible
polynomial.

Then $\cA$ is Faddeev equivalent to a cyclic $F$-algebra of degree
$n$.
\end{prop}

\begin{remark} Cases (i) and (ii) are none other than Propositions
2.2 and 2.3 of \cite{KRTY} (stated in a slightly different form);
for the sake of completeness and reader's convenience we reproduce
the proofs in Section \ref{prop-proof} below. We did not manage to
treat case (iii) in \cite{KRTY} and fill this gap here (for the
case $n=2$ this was done in \cite{RST}).
\end{remark}

We then continue the proof considering, case by case, all possible
configurations of irreducible components of the ramification
divisor $C$:

(i) $C$ is a union of four lines;

(ii) $C$ is a union of two lines and an irreducible conic;

(iii) $C$ is a union of two irreducible conics;

(iv) $C$ is a union of a line and an irreducible cubic curve;

(v) $C$ is irreducible.

In each case we show that by a projective transformation of the
plane one can bring a given algebra $\cD$ to the required form:

\begin{prop} \label{prop:quart}
Let $k$ be an algebraically closed field of characteristic zero,
$\cA$ an algebra over $k({\mathbb P}^2)=k(x,y)$ with ramification
at a plane quartic $C$. Then, after a projective transformation of
${\mathbb P}^2$, $\cA$ viewed as an algebra over $L=l(y)$, where
$l=k(x)$, satisfies the hypotheses of Proposition
$\ref{prop:key}$.
\end{prop}

{\noindent {\it Proof of Theorem $\ref{th:main}$ from Propositions
$\ref{prop:key}$ and $\ref{prop:quart}$}. $\cD$ is Faddeev
equivalent and thus, by Remark \ref{rem:Tsen}, Brauer equivalent
to a cyclic algebra $\cC$ of degree $n$. Hence $\ind (\cC )|n$.
But $n=\exp (\cC )|\ind (\cC )$, implying $\ind (\cC )=n=\deg (\cC
)$. Hence $\cC$ is a division algebra, so $\cD \approx \cC$ and
$\ind (\cD )=\exp (\cD )$. \qed

\section{Proof of Proposition \ref{prop:key}} \label{prop-proof}

Recall that if two algebras have the same ramification at all
finite points, then they have the same ramification at infinity,
by Faddeev's reciprocity law, and thus are Faddeev equivalent.

(i) Using the automorphism group of $\PP^1_K$, we may assume that
the points of ramification of $\cA$ are $t$, $t-1$, and $\infty$
with ramification $K(\sqrt[n]{b_1})/K$, $K(\sqrt[n]{b_2})/K$, and
$K(\sqrt[n]{b_3})/K$ respectively.

The cyclic algebra $\cB = (b_2 t, -b_1^{-1}(t-1))$ has prescribed
ramification at $t$ and $t-1$. Hence $\cA\eqf\cB$.

(ii) Let the ramification of $\cA$ at a quadratic polynomial $f$
be given by a cyclic extension
$K(\theta)(\sqrt[n]{u+v\theta})/K({\theta})$, where $\theta$ is a
root of $f$, $u,v\in K$.

First suppose $v=0$. Consider the algebra $\cB=(u,f)_n$. It is
ramified at $f$ and perhaps at infinity. Its ramification at $f$
coincides with the ramification of $\cA$ at $f$. Hence
$\cA\eqf\cB$.

Now suppose $v\neq 0$. Denote by $c_{uv}$ the value of $f(t)$ at
$-u/v$. We prove that $\cA$ is Faddeev equivalent to $\cB =
(u+vt,f(t)/c_{uv})$. Indeed, $\cB$ could only ramify at $f$,
$t+u/v$, and infinity, but we see that $\cB$ is unramified at
$t+u/v$ since $f(-u/v)/c_{uv}=1$. It has the same ramification at
$f$ as $\cA$ and thus $\cA\eqf\cB$.

(iii) Assume $\cA$ has  ramification $K(\sqrt[n]{a})/K$ at $f$,
where $a\in K(\theta)$, $\theta$ is a root of $f$. Then
$a=a_0+a_1\theta+a_2\theta^2$. We shall consider separately three
cases.

1. $a_1=a_2=0$. Then the algebra $(f(t),a_0)$ has at $f$ the same
ramification as $\cA$, so is Faddeev equivalent to $\cA$.

2. $a_2=0$, $a_1\ne0$. Consider the algebra
$$
\cB=(f(t)/f(-a_0/a_1), a_0+a_1t).
$$
The algebra $\cB$ has at $f$ the same ramification as $\cA$.
Moreover, as in (ii), $\cB$ has no ramification at $a_0+a_1t$.
Hence $\cA\eqf\cB$.

3. $a_2\ne0$. First, let us show that there exists $c\in K$ such
that $a=(b_0+b_1\theta)/(c+\theta)$ for some $b_0,b_1\in K$.

Let $f(t)=t^3-m_2t^2-m_1t-m_0$. We have
$a(c+\theta)=a_2\theta^3+(a_2c+a_1)\theta^2+$ terms of lower
degree in $\theta$. But $\theta^3=m_2\theta^2+m_1\theta +m_0$, so
setting $c=-(m_2+a_1/a_2)$ yields the $\theta^2$-term
$=a_2m_2\theta^2+a_2c\theta+a_1\theta^2=0$, whence
$a(c+\theta)=b_0+b_1\theta$ and $a=(b_0+b_1\theta)/(c+\theta)$.

Consider two subcases:

a) Let $b_1\ne 0$. Consider the algebra
$$
\cC_1=((c+t)^{2n-3}f(t))/((c-b_0/b_1)^{2n-3}f(-b_0/b_1)),
(b_0+b_1t)(c+t)^{n-1}).
$$
Then the algebra $\cC_1$ has at $f$ the same ramification as $\cA$
since
$$
(b_0+b_1\theta)(c+\theta)^{n-1}\equiv (b_0+b_1\theta )(c+\theta
)\min (c+\theta )^n\equiv a(c+\theta )^n \equiv a \pmod
{(K(\theta)^*)^n}.
$$
Furthermore, since both sides of the symbol algebra $\cC_1$ are
polynomials of degree divisible by $n$, using the tame symbol
formula we conclude that $\cC_1$ has no ramification at infinity.
Moreover, $\cC_1$ has no ramification at $t+b_0/b_1$. By the
reciprocity law, $\cC_1$ has no ramification at $c+t$. Thus
$\cA\eqf\cC_1$.

b) If $b_1= 0$, then $\cA\eqf (f(t)/f(-c), b_0(c+t)^{n-1}) $. \qed

\begin{remark}
In each case in the proof of Proposition \ref{prop:key}, we have
provided an explicit cyclic algebra Faddeev equivalent to $\cA$.
\end{remark}


\section{Proof of Proposition \ref{prop:quart}} \label{sec:quart}


(i) Suppose $C$ is the union of 4 lines. Using the automorphism
group of $\PP_k^2$, we may throughout assume that $\cA$ ramifies
at the infinite line and three finite lines. Moreover, we may and
shall assume that one finite line is defined by the equation
$x=0$. Let two other finite lines be defined by the equations
$$
a_ix+b_iy+c_i=0,\;i=1,2.
$$
Hence the algebra $\cA$  viewed as an algebra over $k(x)({\mathbb
P}_{k(x)}^1)$, ramifies at most at two linear polynomials and
perhaps at the infinite point and thus satisfies the hypotheses of
Proposition \ref{prop:key}(i).

(ii) Suppose $C$ is the union of two lines and an irreducible
conic. Using the automorphism group of $\PP_k^2$ as in (i), we may
assume that $\cA$ ramifies at the infinite line, one finite line,
and an irreducible conic. Moreover, we may and shall assume that
the finite line is defined by the equation $x=0$. Let the
conic be defined by the equation $F(x,y)=0$, where $\deg
F(x,y)=2$. Hence the algebra $\cA$  viewed as an algebra over
$k(x)({\mathbb P}_{k(x)}^1)$, ramifies at most at a quadratic
polynomial and perhaps at the infinite point and thus
satisfies the hypotheses of Proposition \ref{prop:key}(ii).

(iii) Again, we may assume that $\cA$ ramifies at the infinite
line and an irreducible cubic defined by the equation
\begin{equation} \label{equat:cubic}
a_0y^3+a_1y^2x+a_2yx^2+a_3x^3+g(x,y)=0,
\end{equation}
where $a_i\in k$ and $\deg\; g(x,y)\le 2$.

Let $a\in k$ be such that $a_0+a_1a+a_2a^2+a_3a^3=0$. Then after
the substitution $x=x'+ay$ we may assume that the cubic defined by
equation (\ref{equat:cubic}) satisfies $a_0=0$.

Hence the algebra $\cA$  viewed as an algebra over $k(x)({\mathbb
P}_{k(x)}^1)$, ramifies at most at a quadratic polynomial in $y$
and perhaps at the infinite point, and we are done.

(iv) Suppose $C$ is the union of two irreducible conics defined by
the homogeneous equations

\begin{equation} \label{equat:2conics1}
a_0y^2+a_1yx+a_2x^2+a_3z^2+a_4xz+a_5yz=0,
\end{equation}
\begin{equation} \label{equat:2conics2}
b_0y^2+b_1yx+b_2x^2+b_3z^2+b_4xz+b_5yz=0,
\end{equation}
where $a_i,b_i\in k$. Since any two curves in the projective plane
intersect, there exists a nontrivial solution
$(u_{11},u_{12},u_{13})$ of the following system of equations
$$
\left\{
\begin{array}{ll}
a_0u_{22}^2+a_1u_{12}u_{22}+a_2u_{12}^2+a_3u_{32}^2+
a_4u_{12}u_{32}+a_5u_{22}u_{32}=0,\\
b_0u_{22}^2+b_1u_{12}u_{22}+b_2u_{12}^2+b_3u_{32}^2+
b_4u_{12}u_{32}+b_5u_{22}u_{32}=0.
\end{array}
\right.
$$
Hence there exist $u_{ij}\in k$, $i=2,3, j=1,2,3$, such that the
matrix $M=(u_{ij})_{i,j=1,2,3}$ is invertible. Applying the
automorphism of $\PP_k^2$ defined by $M^{-1}$, we obtain that the
conics are defined by equations (\ref{equat:2conics1}) and
(\ref{equat:2conics2}) with $a_0=b_0=0$.

Hence the algebra $\cA$, viewed as an algebra over $k(x)({\mathbb
P}_{k(x)}^1)$, ramifies at most at two linear polynomials and
perhaps at the infinite point and, as in (ii), satisfies the hypotheses
of Proposition \ref{prop:key}(ii).


(v) Suppose $C$ is an irreducible quartic.

We need some more notation. Let $X$ be a projective irreducible
curve (not necessarily smooth), $P$ a point on $X$,  $R=O_{X,P}$
the local ring of $X$ at $P$. For $f\in R$ one can define the order
function $\ord_P(f)$ as $\dim_kR/fR$ (see \cite[1.2]{Fu}), continue
this function to $k(X)^*$, and show that $\ord_P\colon k(X)^*\to \Z$
is a homomorphism of abelian groups (see \cite[A.3]{Fu}).

\begin{lemma}  \label{lem:n-th}
Let $C\subset \PP^2_k$ be a projective curve. Suppose that an affine
part $C'=C\cap \A^2_k$ of $C$ is given by an irreducible polynomial
$F(x,y)$ such that $\deg_yF>0$ and $F$ is monic as a polynomial in $y$.
Let $f\in k(C)^*$ be such that $n|\ord_P(f)$ for any $P\in C'$.
Then $N_{k(C)/k(x)}(f)\in (k(x)^*)^n$.
\end{lemma}
\noindent {\it {Proof}}. Consider the projection $\psi\colon
C'\to \A^1_{k(x)}$ corresponding to the homomorphism
$k[x]\to k[x,y]\to k[C']=k[x,y]/(F)$. By hypothesis,
$F$ is monic as a polynomial in $y$, hence the image of $y$ in
$k[C']$ is integral over
$k[x]$. Therefore the ring $k[C']$ is integral over $k[x]$ and thus
$\psi$ is a finite morphism.
Hence $\psi$ is proper and surjective, and by \cite[Prop.~1.4(b)]{Fu} we have
\begin{equation} \label{eq:Fu}
\psi_*(\div (f)) = \div (N_{k(C)/k(x)}(f)).
\end{equation}
By hypothesis, all coefficients of $\div (f)$ at finite points $P$
are divisible by $n$. We conclude that the coefficients of the divisor
in the right-hand side of (\ref{eq:Fu}) at all finite points are also
divisible by $n$. This means that the values of the corresponding
valuations are divisible by $n$, and since $k$ is algebraically
closed, $N_{k(C)/k(x)}(f)$ is an $n$-th power. \qed

\bigskip

We need another auxiliary result.

\begin{lemma} \label{prop:OneCurve}
Let $k$ be an algebraically closed field of characteristic zero,
$l=k(x)$, $L=k(\PP^2_k)=k(x,y)$,
and $\cA$ an algebra over $L$ with ramification at
only one irreducible curve $C\subset \PP^2$
satisfying the hypotheses of Lemma $\ref{lem:n-th}.$
Then $\cA$ as an
algebra over $l(\PP^1_l)=k(x)(y)$ has no ramification at the
infinite point.
\end{lemma}

{\noindent{\it Proof.} Maintain the notation of Lemma \ref{lem:n-th}.
Denote $n=\exp (\cA)$. Suppose that the
ramification of $\cA$ at $C$ is defined by $f\in k(C)$. Note that
$k(C)$ is the residue field corresponding to the valuation of
$k(x,y)$ defined by $C$. According to \cite{AM} (see also
\cite[p.196]{Sa95} for more general setting), there is an exact
sequence
\begin{equation} \label{eq:ram12}
0 \rightarrow \Br(L) \stackrel{\oplus\ram_Y}{\longrightarrow}
\bigoplus_{Y\subset \A_k^{2(1)}} H^1(k(Y),\Q/\Z)
\stackrel{\phi}{\longrightarrow} \bigoplus_{Q\in \A_k^{2(2)}}
\Q/\Z \longrightarrow 0,
\end{equation}
where $Y$ and $Q$ range over all irreducible subvarieties of
$\A_k^2$ of codimension 1 and 2, respectively. The map $\phi$ is
defined as follows. Let $Q\in \A_k^{2(2)}$, $Y\in \A_k^{2(1)}$.
Suppose $\alpha \in H^1(k(Y),\Q/\Z)$ is of order $n$ and defined
by $a\in k(Y)$. Assume $Q\in Y$. Set $r=\ord_Q(a)$ and define the
ramification of $\alpha$ at $Q$ by $\ram_Q(\alpha)=r/n$.

If $Q$ is not a subset of $Y$, we set $\ram_Q(\alpha)=0$. The map
$\phi$ is the sum of all maps $\ram_Q$.

Using sequence (\ref{eq:ram12}), we obtain
$\ram_P(\ram_C([\cA]))=0$ for any $P\in C'$. Hence
$n|\ord_P(f)$ for any $P\in C'$. By Lemma
\ref{lem:n-th}, we have $N_{k(C)/k(x)}(f)\in (k(x)^*)^n$. Hence
from the exact sequence (\ref{eq:ramFad}) we obtain $\cor_p(f)=0$,
where $p$ is the closed point of $\PP_{k(x)}^1$ corresponding to
the valuation of $k(x,y)$ defined by $C$. Thus $\cA$ as an algebra
over $l(\PP^1_l)=k(x)(y)$ has no ramification at the infinite
point. \qed

\bigskip

We can now prove the assertion of case (v) of Proposition
\ref{prop:quart}.

Without loss of generality we may assume that
$\deg_yF>0$ for otherwise $\cA$ is unramified over
$k(x)(\PP^1_{k(x)})$. Moreover, we may assume
$\deg_yF>2$ for otherwise we are reduced to Proposition
\ref{prop:key}((i) or (ii)). Furthermore, by an argument analogous to that
in case (iii), we may and shall assume that the coefficient at $y^4$ is zero.
Let us show that one can transform $C$ so that it satisfies the hypotheses
of Lemma \ref{lem:n-th}.

Write an affine equation for $C'$ in the form
$$
y^3(a_1x +a_0) + y^2g_2(x) + y g_3(x) + g_4(x) = 0,
$$
where $g_i(x)\in k[x]$ denotes a polynomial of degree $i$ ($2\leq i\leq 4$).
If $a_1=0$, on dividing by $a_0$ we arrive at the hypotheses of Lemma
\ref{lem:n-th}. So suppose $a_1\neq 0$ and consider the homogeneous
equation of $C$:
$$
y^3(a_1x+a_0z) + y^2G_2(x,z) +y G_3(x,z) + G_4(x,z) =0,
$$
where $G_i(x,z)$ stands for the homogenization of $g_i(x)$.
The projective transformation of the plane
$x'=a_1x+a_0z$, $y=y$, $z=z$ brings $C$ to the form
$$
y^3x' +y^2H_2(x',z) + yH_3(x',z) + H_4(x',z) =0,
$$
where $H_i$ is homogeneous of degree $i$. Going over to another affine
chart on dividing by ${x'}^4$, we obtain an equation for another affine
part of $C$ which satisfies the hypotheses of Lemma \ref{lem:n-th}.

We can now apply Lemma \ref{prop:OneCurve} to conclude that the algebra
$\cA$  viewed as an
algebra over $k(x)({\mathbb P}_{k(x)}^1)$, ramifies only at a cubic
irreducible polynomial, the hypothesis of Proposition
\ref{prop:key}(iii).

Proposition \ref{prop:quart} is proved. This completes the proof
of Theorem \ref{th:main}. \qed





\enddocument